# Chaos predictability in a chemical reactor

Marek Berezowski

*Cracow University of Technology ,Poland*

**Abstract**

The dynamics of the tubular chemical reactor with mass recycle were examined. In such a system, temperature and concentrations may oscillate chaotically. This means that state variable values are then unpredictable. In this paper it has been shown that despite the chaos, the behaviour of such a reactor can be predictable. It has been shown that this phenomenon can occur in two cases. The first case concerns intermittent chaos. It has been shown that intermittent outbursts can occur at regular intervals. The second case concerns transient chaos, i.e. a situation when the chaos occurs only for a certain period of time, e.g. only during start-up. This phenomenon makes it impossible to predict what will occur in the reactor in the nearest time, but, makes it possible to precisely determine the values of the variables even in the distant future. Both of these phenomena were tested by numerical simulation of the mathematical model of the reactor.

## 1. Introduction

The unpredictability of phenomena is a fundamental feature of chaos. This means that even the smallest deviation from the initial conditions causes, after a certain period of time - called the Lyapunov time, a significant change in the trajectory position in the phase space [Berezowski, 2009]. In physical and industrial processes this is manifested by the fact that two time series evolve quite differently, even though they were practically identical in the beginning. In some cases – as discussed in the paper, the principle of unpredictability may be infringed, even though the changes are chaotic.

These phenomena are shown on the example of analyzing a model of a non-adiabatic tubular chemical reactor with mass recycle (Fig. 1). This model has already been tested many times in the scientific literature, where it has been shown that temperature-concentrations oscillations can be chaotic, which is mainly due to recirculation coupling [for example: Jacobsen & Berezowski,1998].

One of the types of chaos (which also occurs in this reactor) is the so called: "intermittency" [Pomeau & Manneville, 1980]. Figuratively speaking, it means that "regular" changes, preceded by rapid bursts, occur in the time series. So far, the theory and practice have proved that the moments in which the bursts occur, are unpredictable.

The scope of the paper is to indicate that in some cases, the above statement may not always hold. The analysis of a chemical reactor model reveals that under certain conditions intermittent bursts may occur (practically) at regular time intervals. This phenomenon does not breach the theory of chaos, because the values of the amplitudes of the bursts, as well as their duration, remain unpredictable. This phenomenon has been shown to occur at the boundary between the chaotic and periodic solution. Knowledge of this phenomenon is very valuable in the industrial process because it allows to determine the moments of time when explosions occur (e.g. high temperatures) and thus to avoid their harmful effects.

The second interesting phenomenon shown in this paper is the anticipation of changes in the so-called transient chaos. It occurs when chaos occurs only at specific time intervals, and not throughout the entire time series. For example, if chaos occurs only during start-up, this phenomenon makes it impossible to predict what will occur in the reactor in the nearest time, but, makes it possible to precisely determine the values of the variables even in the distant future.

## 2. The model of a chemical reactor

A mathematical model of a non-adiabatic homogeneous tubular chemical reactor without dispersion and with outlet flow recirculation was considered in the paper. (Fig. 1) [Jacobsen & Berezowski, 1998], [Berezowski, 2000], [Berezowski, 2001], [Berezowski, 2003], [Berezowski & Bizon, 2006], [Berezowski, 2006], [Berezowski, 2009], [Berezowski, 2013], [Luss & Amundson, 1967], [Reilly & Schmitz, 1966], [Reilly & Schmitz, 1967], [Lauschke & Gilles, 1994], [Femat et al., 2004], [Kienle et al., 1995], [Antoniades & Christofides, 2001].

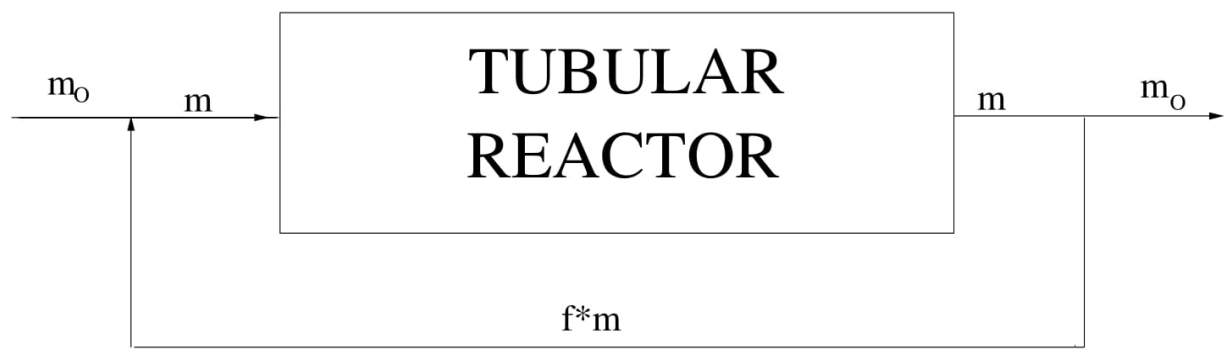

Fig. 1. Schematic diagram of the reactor with recycle

The corresponding balance equations are as follows:

mass balance:

$$\frac{\partial \alpha}{\partial \tau} + \frac{\partial \alpha}{\partial \varsigma} = \phi(\alpha, \Theta) \qquad (1)$$

heat balance:

$$\frac{\partial \Theta}{\partial \tau} + \frac{\partial \Theta}{\partial \varsigma} = \phi(\alpha, \Theta) + (1-f)\delta(\Theta_H - \Theta) \qquad (2)$$

where: $0 \leq \varsigma \leq 1$ [Jacobsen & Berezowski, 1998; 2000; 2001; 2003; 2006; 2009, 2013; 2016], [Berezowski & Bizon, 2006].

Assuming that a single reaction of the $A \rightarrow B$ *n-th* order occurs in the reactor, the kinetic function of the reaction has the form:

$$\phi(\alpha, \Theta) = (1-f)\text{Da}(1-\alpha)(1-\alpha)^n \exp(\gamma \frac{\beta\Theta}{1+\beta\Theta}). \qquad (3)$$

Due to the presence of the recirculation loop, the mathematical model of the above mentioned reactor must be supplemented with appropriate boundary conditions concerning mass and heat:

$$\alpha(\tau, 0) = f\alpha(\tau, 1) \qquad (4)$$

$$\Theta(\tau, 0) = f\Theta(\tau, 1). \qquad (5)$$

From the mathematical point of view, such model is discrete, rendering in its solution a rectangular wave of the changes of state α and Θ. According to available scientific publications, in some cases the wave may have a chaotic character [Jacobsen & Berezowski, 1998], [Berezowski, 2000], [Berezowski, 2001], [Elnashaie & Abashar, 1995].

## 3. Results of computations and their analysis

All calculations were carried out by the method of numerical simulation of the model (1) - (5). To ensure the reliability of the results, various integration methods were used in the calculations, including the Runge-Kutta, Simpson and Euler methods. In each case identical results were obtained. Due to the plug-flow of the stream through the reactor and due to the nature of the boundary conditions (4) - (5), the analyzed model has - from a mathematical point of view - a discrete character. This means that we can observe the state variables at discrete moments of time 0, 1, 2, etc. The time interval = 1 corresponds in this case to the time of the single flow of the sample particle from the inlet to the outlet of the reactor tube. This is called residence time in the reactor. This mechanism has been described in detail, including at work [Jacobsen & Berezowski, 1998].

The computations were carried out for the following values of the reactor parameters: *Da=0.15, n=1.5, β=2, γ=15, δ=3, f=0.5*. The Feigenbaum diagram was designated for these values, as a function of the dimensionless temperature of the cooling medium $\Theta_H$ (Fig. 2) [Feigenbaum , 1980].

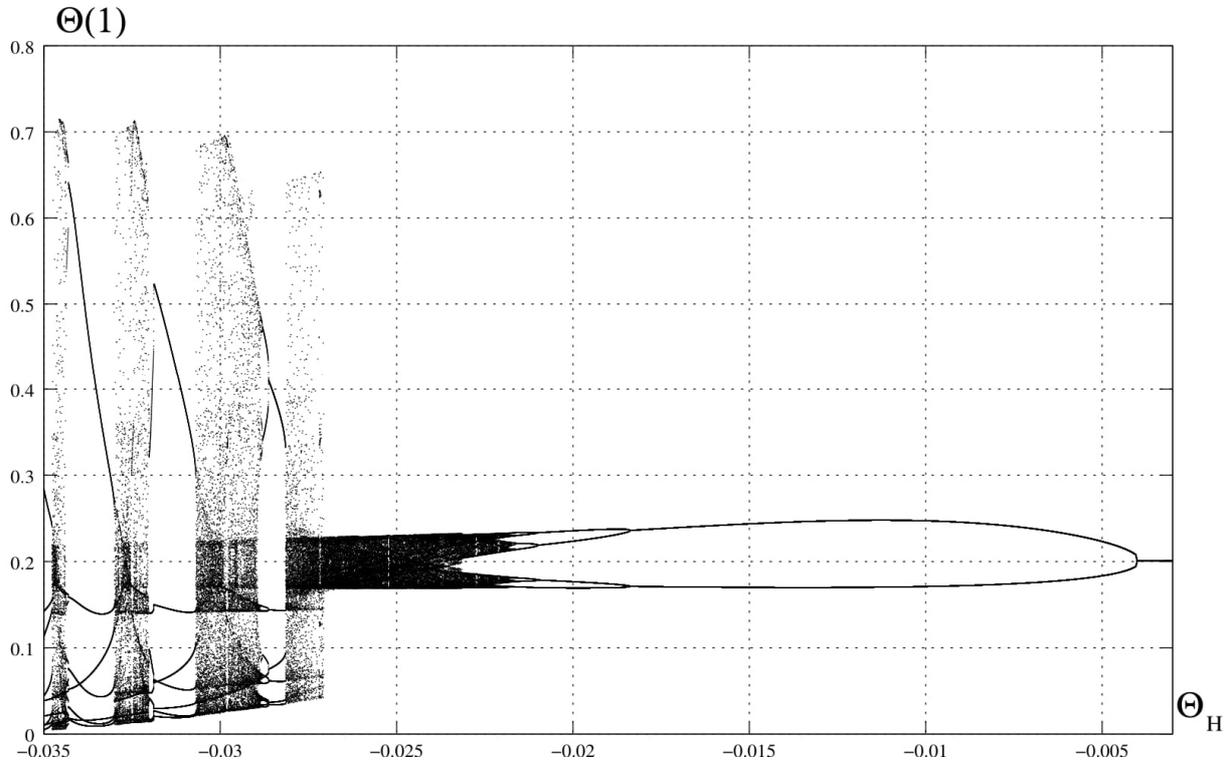

Fig. 2. Feigenbaum diagram

In the analysis of the diagram from the right to the left side it may be clearly seen that, in the initial phase, the scenario of reaching chaos involves the period doubling, where oscillations of relatively low amplitude occur. Going further, a crisis generating oscillations of high amplitude are discerned. This is the second type of crisis, an interior crisis, where the size of the chaotic attractor suddenly increases. This fragment is studied in detail in the following analysis [Berezowski, 2013], [Grebogi et al., 1983].

It may be easily noticed that in the zone of the crisis, windows of chaotic and periodic solutions occur alternatingly, which is a characteristic feature of such type of diagram. It may also be demonstrated that the number of the periodic windows is infinite. In the majority, however, they are too narrow to be visible in the diagram.

Assuming, for example, $\Theta_H$=-0.0335, a periodic chaotic time series is obtained (Fig. 3).

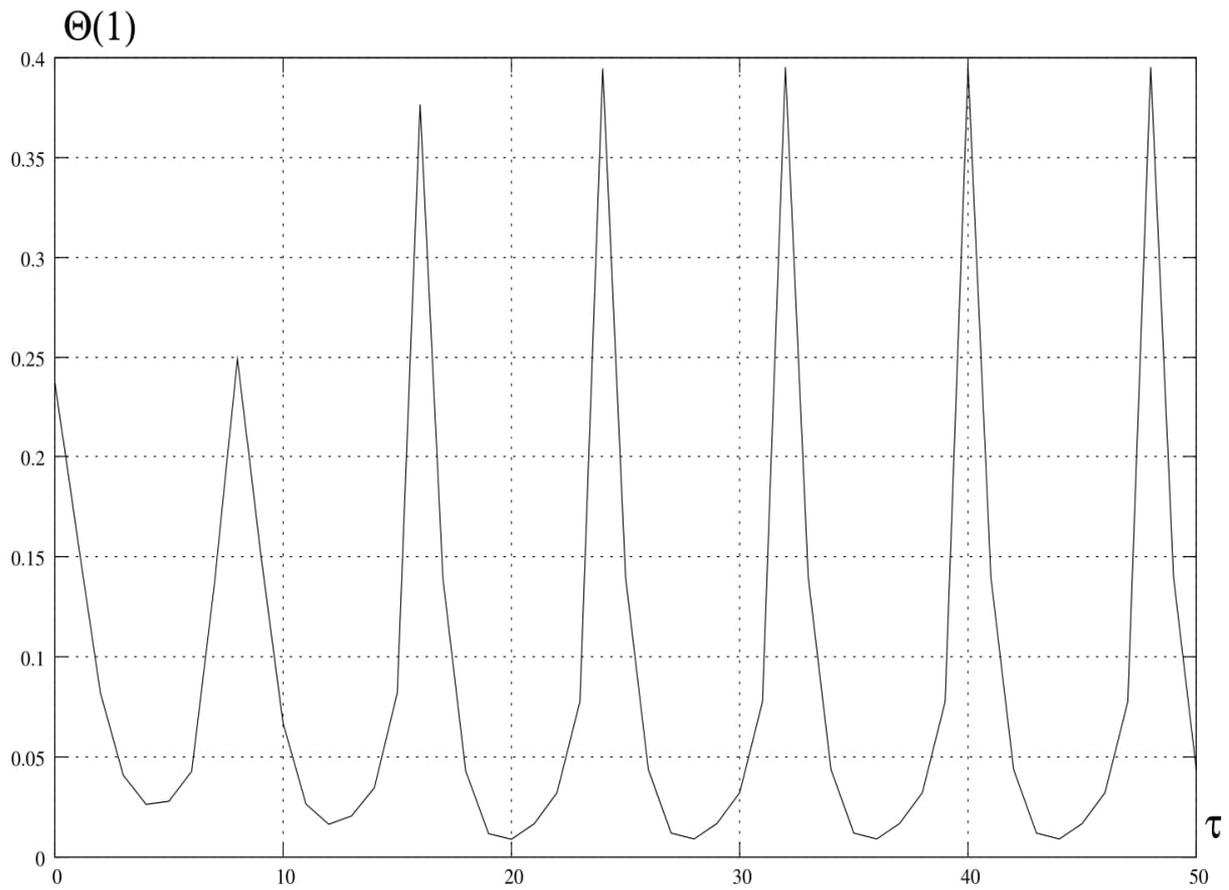

Fig. 3. Time series

It is insensitive to the initial conditions and has a negative Lyapunov exponent. In turn, assuming that $\Theta_H=-0.03299$, a chaotic time series is obtained (Fig. 4). It is sensitive to the initial conditions and has a positive Lyapunov exponent [Jacobsen & Berezowski, 1998], [Berezowski, 2000], [Benettin et al., 1980].

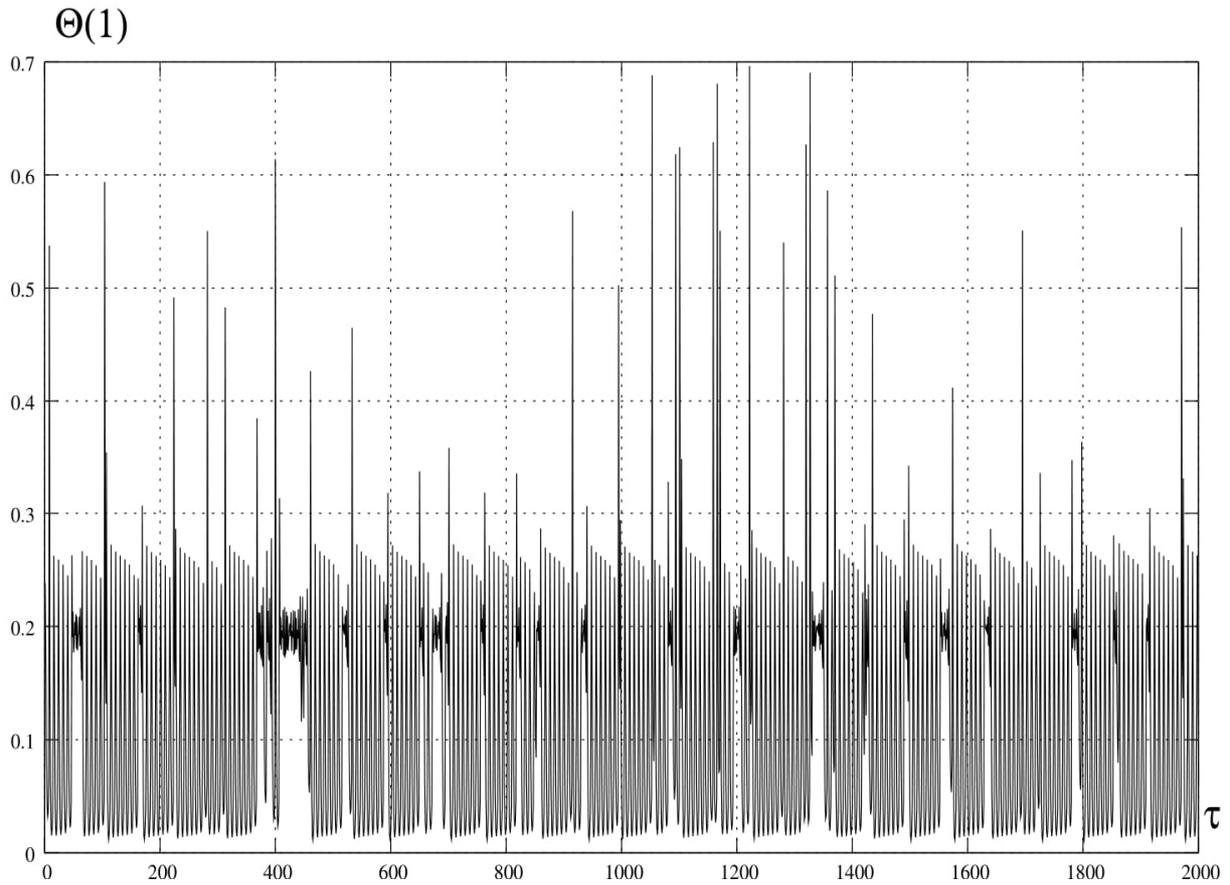

Fig. 4. Non-periodic intermittency

Concurrently, it should be noticed that the series is of an intermittent type, as intervals of regular changes divided by specific bursts are visible. The intervals also occur at uneven time distances, which means that they are impossible to predict.

An interesting phenomenon discussed further in the paper, occurs at the boundary between the chaotic and periodic solution (for example, for $\Theta_H = -0.033003$ - Fig. 5).

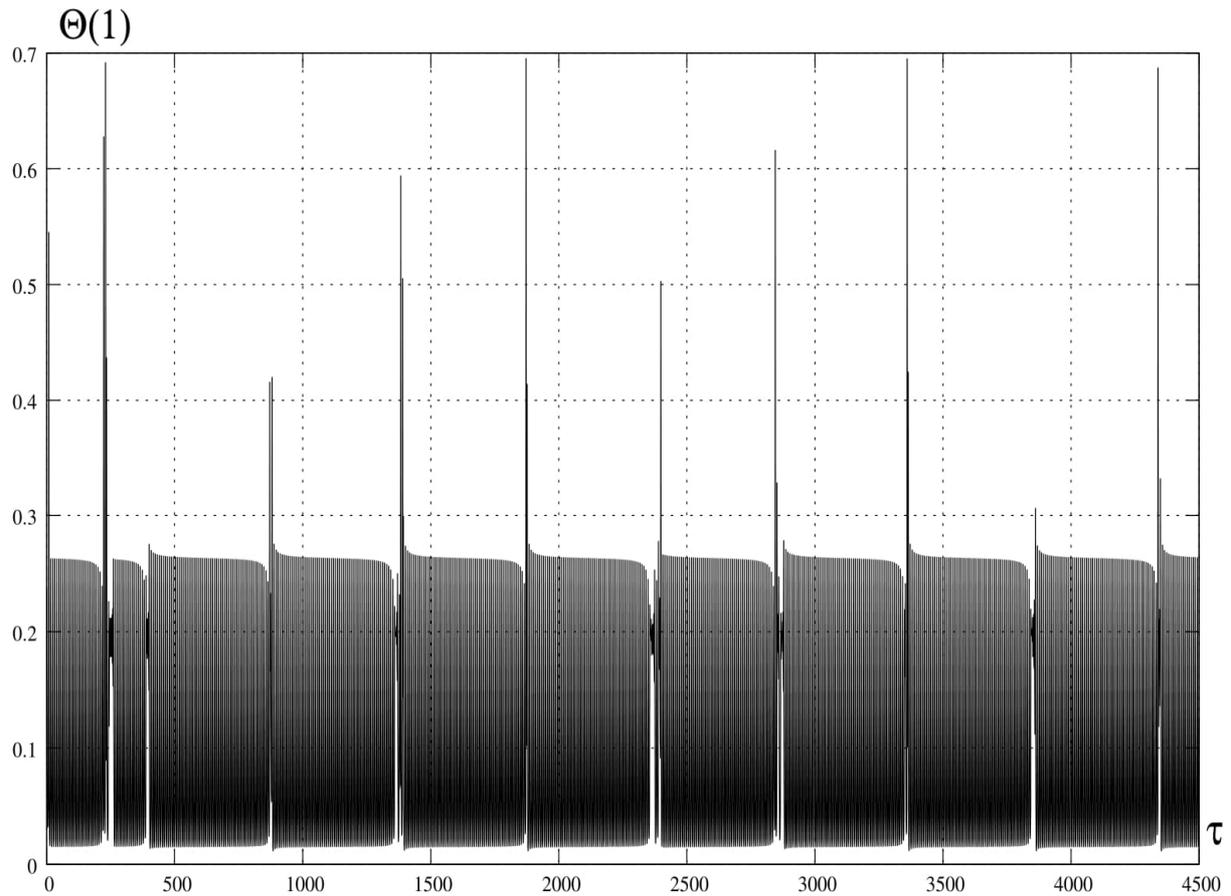

Fig. 5. Short-period periodic intermittency

Intermittency may also be observed there, but, unlike in the previous example – the bursts occur (practically) at regular time intervals, which, in this case, is at every 479 units. This phenomenon renders an opportunity of predicting the time moments at which the successive bursts occur. However, the value of the Lyapunov exponent in this case is positive and equals $\lambda=0.02174$. The amplitudes of the bursts and their duration remain unpredictable.

The Poincaré map presented in Fig.6 is related to intermittent bursts. It may be easily noticed that the map has the Henon attractor character, typical of chaotic discrete systems.

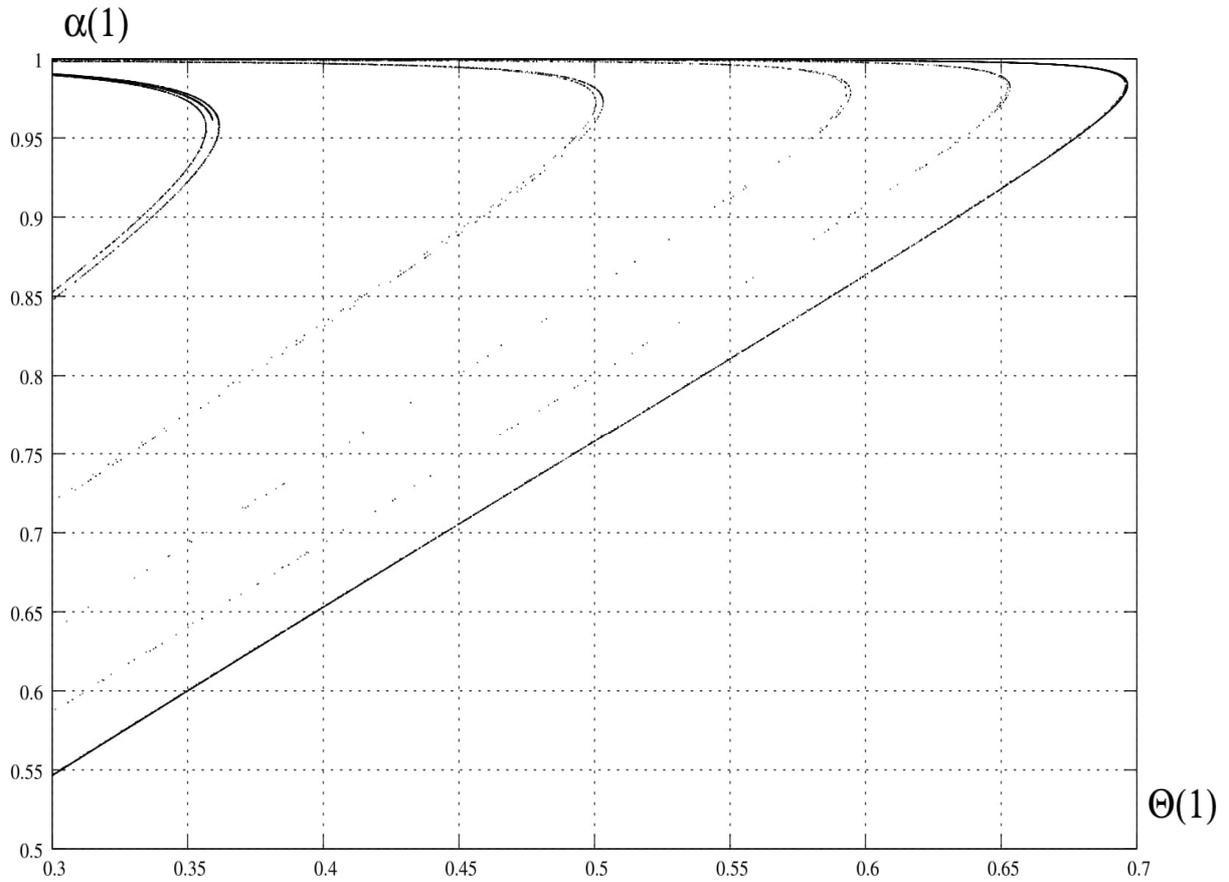

Fig. 6. Poincaré map of intermittent bursts

It should also be added that the phenomenon described above occurs everywhere where chaos osculates directly with periodic oscillations, i.e. at various points of the Feigenbaum diagram in the crisis zone. Using the fragment of the diagram (Fig. 7), another example of regular intermittency is observed, this time for $\Theta_H$=-0.0327049346 (Fig. 8). In such case, the interval between the bursts is longer and amounts to about 3980 units.

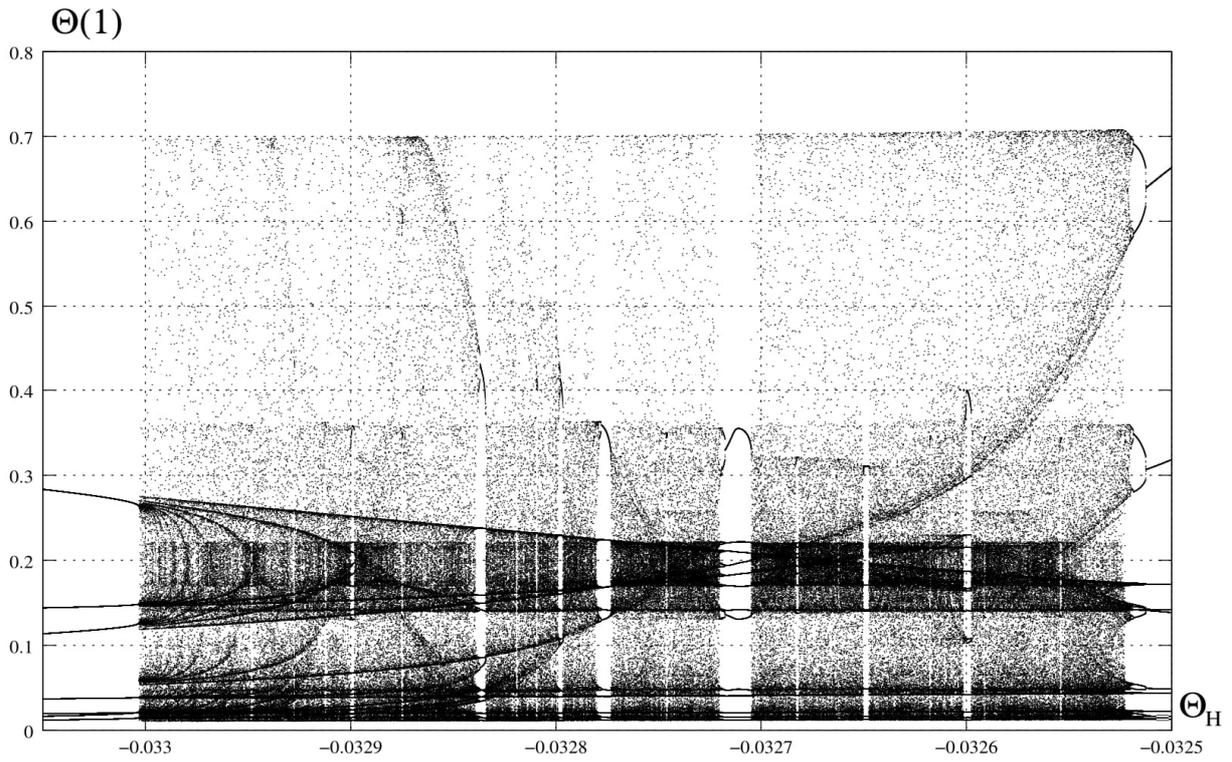

Fig. 7. Fragment of the Feigenbaum diagram

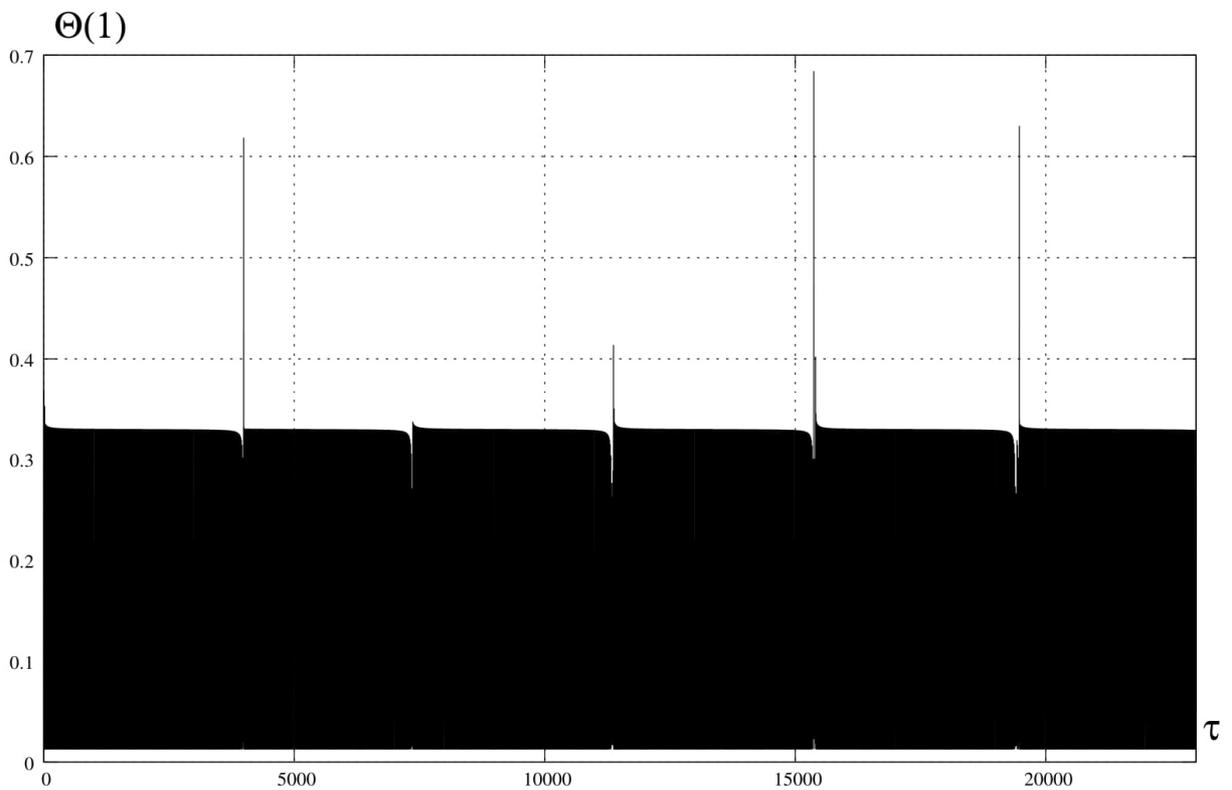

Fig. 8. Long-period periodic intermittency

In addition, another type of predictable chaos is discussed in the paper, for example, for $\Theta_H=-0.0330006$, the so called "transient" chaos occurs (Fig. 9) [Berezowski, 2016]. An interesting feature of this chaos is that it appears in the initial phase of the time series, but decays in the final phase, becoming a periodic oscillations.

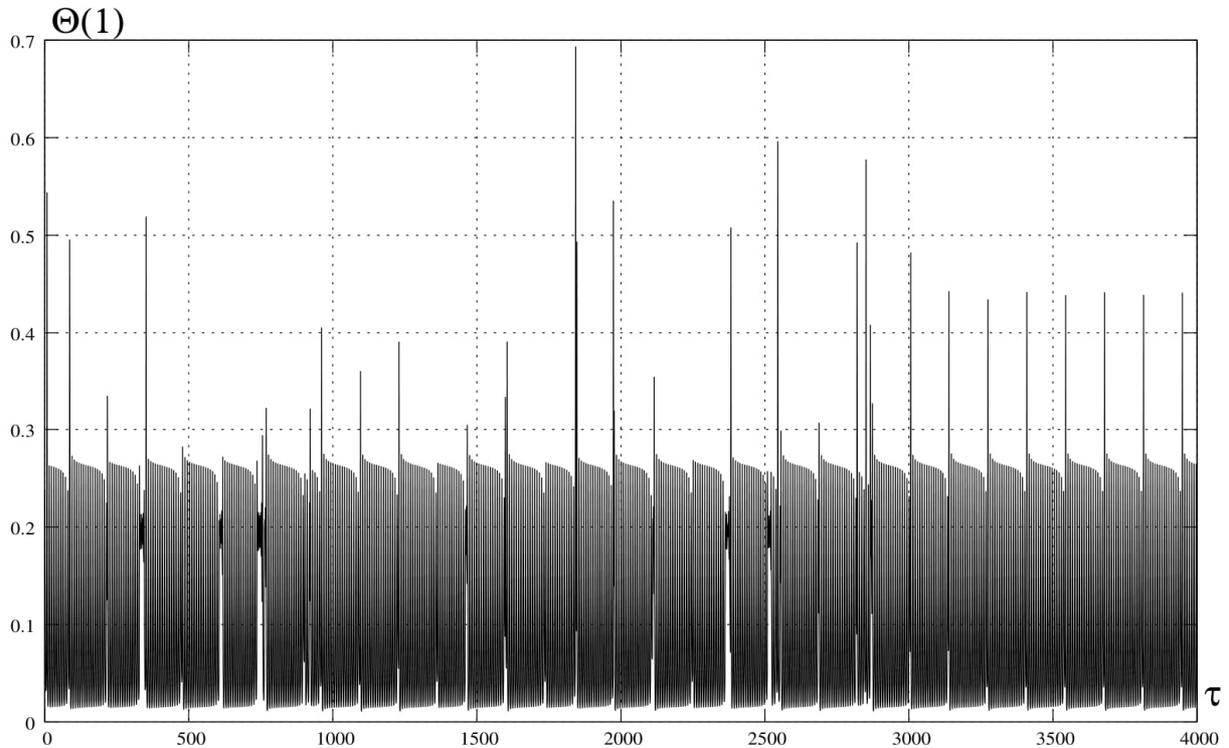

Fig. 9. Transient chaos

Such phenomena occur in physical and industrial processes due to natural causes, for example, as an effect of the catalyst deactivation, or natural decay of a certain parameter in time. Mathematical models of such systems are usually non-autonomous and express time in an explicit form, leading to the degradation of some parameters. This is not the case of the reactor analysed in the paper. The discussed model is autonomous and none of its parameters decay in the passage of time. Despite this, in the analysed example, the time series of the dimensionless temperature $\Theta(1)$ is chaotic only for $\tau<3300$. Then, the series becomes periodic. The chaotic character of the initial phase of the series is also indicated by the sensitivity to the initial conditions and positive value of the Lyapunov exponent: $\lambda=0.063$. The same exponent calculated for very long time approaches zero. However, calculated for $\tau>3300$, it has the value of $\lambda=-0.0006636$, which confirms the periodicity of the final phase. Thus, this phenomenon makes it impossible to predict what will occur in the reactor in the nearest time, but, makes it possible to precisely determine the values of the variables even in the distant future.

## 4. Summary and conclusions

A theoretical analysis of the dynamics of a non-diabatic homogeneous tubular reactor without dispersion and with the outlet mass circulation was conducted in the paper. It was proved that the reactor may generate chaotic oscillations, in which intermittent bursts may occur (practically) at regular time intervals. Hence, the moments at which the bursts appear are predictable. However, the amplitudes of the bursts and their duration remain unpredictable. Furthermore, it should be noticed that their duration is relatively short (one or several impulses). It was also indicated that this phenomenon occurs at the boundary between the chaos and the periodic oscillations. The knowledge of this phenomenon is very important from the process point of view, both for the time of the apparatus operation and for the stage of its design, as it enables the determination of the moments of bursts and, in consequence, possible protection of the system by appropriate control.

Another finding discussed in the paper is the possibility of the existence of another type of predictable chaos, the so called "transient" chaos [Berezowski, 2016]. This phenomenon involves chaotic, therefore unpredictable, operation of the reactor in the initial phase, but periodic, therefore predicable, operation of the reactor in the next phase. In practice, this means that we are not able to predict changes in the temperature and concentrations in the initial phase, but we can precisely designate the values of these parameters at very distant moments. Thus, the knowledge of this phenomenon is very important to the start-up phase of the reactor. Furthermore, it also facilitates the elimination of the control of the apparatus after its start-up phase.

**Symbols**

Da – Damköhler number

f – recycle coefficient

F - volumetric flow, $m^3/s$

L – length of reactor, *m*

m – mass flow, *kg/s*

n – reaction order

V – volume of reactor, $m^3$

t – time, *s*

x – location along the reactor, *m*

*Greek letters:*

α – conversion degree

β – dimensionless enthalpy

ς – dimensionless location along the reactor ($= \frac{x}{L}$)

δ – dimensionless coefficient of heat exchange

γ – dimensionless activation energy

Θ – dimensionless temperature

τ – dimensionless time ($= \frac{tF}{V}$)

*Indices:*

H – referring to the cooling medium

o – refer to feed